# Maximum likelihood estimation of a log-concave density and its distribution function: Basic properties and uniform consistency

LUTZ DÜMBGEN[1] and KASPAR RUFIBACH[2]

[1]*Institute of Mathematical Statistics and Actuarial Science, University of Bern, Sidlerstrasse 5, CH-3012 Bern, Switzerland. E-mail: duembgen@stat.unibe.ch*
[2]*Abteilung Biostatistik, Institut für Sozial- und Präventivmedizin, Universität Zürich, Hirschengraben 84, CH-8001 Zürich, Switzerland. E-mail: kaspar.rufibach@ifspm.uzh.ch*

We study nonparametric maximum likelihood estimation of a log-concave probability density and its distribution and hazard function. Some general properties of these estimators are derived from two characterizations. It is shown that the rate of convergence with respect to supremum norm on a compact interval for the density and hazard rate estimator is at least $(\log(n)/n)^{1/3}$ and typically $(\log(n)/n)^{2/5}$, whereas the difference between the empirical and estimated distribution function vanishes with rate $o_{\mathrm{p}}(n^{-1/2})$ under certain regularity assumptions.

*Keywords:* adaptivity; bracketing; exponential inequality; gap problem; hazard function; method of caricatures; shape constraints

## 1. Introduction

Two common approaches to nonparametric density estimation are smoothing methods and qualitative constraints. The former approach includes, among others, kernel density estimators, estimators based on discrete wavelets or other series expansions and estimators based on roughness penalization. Good starting points for the vast literature in this field are Silverman (1982, 1986) and Donoho *et al.* (1996). A common feature of all of these methods is that they involve certain tuning parameters, for example, the order of a kernel and the bandwidth. A proper choice of these parameters is far from trivial since optimal values depend on unknown properties of the underlying density $f$. The second approach avoids such problems by imposing qualitative properties on $f$, for example, monotonicity or convexity on certain intervals in the univariate case. Such assumptions are often plausible or even justified rigorously in specific applications.

Density estimation under shape constraints was first considered by Grenander (1956), who found that the nonparametric maximum likelihood estimator (NPMLE) $\hat{f}_n^{\mathrm{mon}}$ of a non-increasing density function $f$ on $[0,\infty)$ is given by the left derivative of the least concave majorant of the empirical cumulative distribution function on $[0,\infty)$. This work was continued by Rao (1969) and Groeneboom (1985, 1988), who established asymptotic distribution theory for $n^{1/3}(f - \hat{f}_n^{\mathrm{mon}})(t)$ at a fixed point $t > 0$ under certain regularity conditions and analyzed the non-Gaussian limit distribution. For various estimation problems involving monotone functions, the typical rate of convergence is $O_{\mathrm{p}}(n^{-1/3})$ pointwise. The rate of convergence with respect to supremum norm is further decelerated by a factor of $\log(n)^{1/3}$ (Jonker and van der Vaart (2001)). For applications of monotone density estimation, consult, for example, Barlow *et al.* (1972) or Robertson *et al.* (1988).

Monotone estimation can be extended to cover unimodal densities. Remember that a density $f$ on the real line is unimodal if there exists a number $M = M(f)$ such that $f$ is non-decreasing on $(-\infty, M]$ and non-increasing on $[M, \infty)$. If the true mode is known a priori, unimodal density estimation boils down to monotone estimation in a straightforward manner, but the situation is different if $M$ is unknown. In that case, the likelihood is unbounded, problems being caused by observations too close to a hypothetical mode. Even if the mode was known, the density estimator is inconsistent at the mode, a phenomenon called "spiking". Several methods were proposed to remedy this problem (see Wegman (1970), Woodroofe and Sun (1993), Meyer and Woodroofe (2004) or Kulikov and Lopuhaä (2006)), but all of them require additional constraints on $f$.

The combination of shape constraints and smoothing was assessed by Eggermont and La-Riccia (2000). To improve the slow rate of convergence of $n^{-1/3}$ in the space $\mathrm{L}_1(\mathbb{R})$ for arbitrary unimodal densities, they derived a Grenander-type estimator by taking the derivative of the least concave majorant of an integrated kernel density estimator rather than the empirical distribution function directly, yielding a rate of convergence of $O_{\mathrm{p}}(n^{-2/5})$.

Estimation of a convex decreasing density on $[0,\infty)$ was pioneered by Anevski (1994, 2003). The problem arose in a study of migrating birds discussed by Hampel (1987). Groeneboom *et al.* (2001) provide a characterization of the estimator, as well as consistency and limiting behavior at a fixed point of positive curvature of the function to be estimated. They found that the estimator must be piecewise linear with knots between the observation points. Under the additional assumption that the true density $f$ is twice continuously differentiable on $[0,\infty)$, they show that the MLE converges at rate $O_{\mathrm{p}}(n^{-2/5})$ pointwise, somewhat better than in the monotone case. Monotonicity and convexity constraints on densities on $[0,\infty)$ have been embedded into the general framework of $k$–monotone densities by Balabdaoui and Wellner (2008). In a technical report, we provide a more thorough discussion of the similarities and differences between $k$-monotone density estimation and the present work (Dümbgen and Rufibach (2008)).

In the present paper, we impose an alternative, and quite natural, shape constraint on the density $f$, namely, log-concavity. That means

$$f(x) = \exp\varphi(x)$$

for some concave function $\varphi:\mathbb{R} \to [-\infty, \infty)$. This class is rather flexible, in that it generalizes many common parametric densities. These include all non-degenerate normal densities, all Gamma densities with shape parameter $\geq 1$, all Weibull densities with exponent $\geq 1$ and all beta densities with parameters $\geq 1$. Further examples are the logistic and Gumbel densities. Log-concave densities are of interest in econometrics; see Bagnoli and Bergstrom (2005) for a summary and further examples. Barlow and Proschan (1975) describe advantageous properties of log-concave densities in reliability theory, while Chang and Walther (2007) use log-concave densities as an ingredient in nonparametric mixture models. In nonparametric Bayesian analysis, too, log-concavity is of certain relevance (Brooks (1998)).

Note that log-concavity of a density implies that it is also unimodal. It will turn out that by imposing log-concavity, one circumvents the spiking problem mentioned before, which yields a new approach to estimating a unimodal, possibly skewed density. Moreover, the log-concave density estimator is fully automatic, in the sense that there is no need to select any bandwidth, kernel function or other tuning parameters. Finally, simulating data from the estimated density is rather easy. All of these properties make the new estimator appealing for use in statistical applications.

Little large sample theory is available for log-concave estimators thus far. Sengupta and Paul (2005) considered testing for log-concavity of distribution functions on a compact interval. Walther (2002) introduced an extension of log-concavity in the context of certain mixture models, but his theory does not cover asymptotic properties of the density estimators themselves. Pal *et al.* (2006) proved the log-concave NPMLE to be consistent, but without rates of convergence.

Concerning the computation of the log-concave NPMLE, Walther (2002) and Pal *et al.* (2006) used a crude version of the iterative convex minorant (ICM) algorithm. A detailed description and comparison of several algorithms can be found in Rufibach (2007), while Dümbgen *et al.* (2007a) describe an active set algorithm, which is similar to the vertex reduction algorithms presented by Groeneboom *et al.* (2008) and seems to be the most efficient one at present. The ICM and active set algorithms are implemented within the R package "logcondens" by Rufibach and Dümbgen (2006), accessible via "CRAN". Corresponding MATLAB code is available from the first author's homepage.

In Section 2, we introduce the log-concave maximum likelihood density estimator, discuss its basic properties and derive two characterizations. In Section 3, we illustrate this estimator with a real data example and explain briefly how to simulate data from the estimated density. Consistency of this density estimator and the corresponding estimator of the distribution function are treated in Section 4. It is shown that the supremum norm between estimated density, $\hat{f}_n$, and true density on compact subsets of the interior of $\{f > 0\}$ converges to zero at rate $O_{\mathrm{p}}((\log(n)/n)^\gamma)$, with $\gamma \in [1/3, 2/5]$ depending on $f$'s smoothness. In particular, our estimator adapts to the unknown smoothness of $f$. Consistency of the density estimator entails consistency of the distribution function estimator. In fact, under additional regularity conditions on $f$, the difference between the empirical c.d.f. and the estimated c.d.f. is of order $o_{\mathrm{p}}(n^{-1/2})$ on compact subsets of the interior of $\{f > 0\}$.

As a by-product of our estimator, note the following. Log-concavity of the density function $f$ also implies that the corresponding hazard function $h = f/(1 - F)$ is non-decreasing (cf. Barlow and Proschan (1975)). Hence, our estimators of $f$ and its c.d.f. $F$ entail a consistent and non-decreasing estimator of $h$, as pointed out at the end of Section 4.

Some auxiliary results, proofs and technical arguments are deferred to the Appendix.

## 2. The estimators and their basic properties

Let $X$ be a random variable with distribution function $F$ and Lebesgue density

$$f(x) = \exp \varphi(x)$$

for some concave function $\varphi : \mathbb{R} \to [-\infty, \infty)$. Our goal is to estimate $f$ based on a random sample of size $n > 1$ from $F$. Let $X_1 < X_2 < \cdots < X_n$ be the corresponding order statistics. For any log-concave probability density $f$ on $\mathbb{R}$, the normalized log-likelihood function at $f$ is given by

$$\int \log f \, d\mathbb{F}_n = \int \varphi \, d\mathbb{F}_n, \tag{1}$$

where $\mathbb{F}_n$ stands for the empirical distribution function of the sample. In order to relax the constraint of $f$ being a probability density and to get a criterion function to maximize over the convex set of *all* concave functions $\varphi$, we employ the standard trick of adding a Lagrange term to (1), leading to the functional

$$\Psi_n(\varphi) := \int \varphi \, d\mathbb{F}_n - \int \exp \varphi(x) \, dx$$

(see Silverman (1982), Theorem 3.1). The nonparametric maximum likelihood estimator of $\varphi = \log f$ is the maximizer of this functional over all concave functions,

$$\hat{\varphi}_n := \underset{\varphi \text{ concave}}{\arg \max} \, \Psi_n(\varphi)$$

and $\hat{f}_n := \exp \hat{\varphi}_n$.

*Existence, uniqueness and shape of $\hat{\varphi}_n$.* One can easily show that $\Psi_n(\varphi) > -\infty$ if and only if $\varphi$ is real-valued on $[X_1, X_n]$. The following theorem was proven independently by Pal *et al.* (2006) and Rufibach (2006). It also follows from more general considerations in Dümbgen *et al.* (2007a), Section 2.

**Theorem 2.1.** *The NPMLE $\hat{\varphi}_n$ exists and is unique. It is linear on all intervals $[X_j, X_{j+1}]$, $1 \leq j < n$. Moreover, $\hat{\varphi}_n = -\infty$ on $\mathbb{R} \setminus [X_1, X_n]$.*

*Characterizations and further properties.* We provide two characterizations of the estimators $\hat{\varphi}_n$, $\hat{f}_n$ and the corresponding distribution function $\hat{F}_n$, that is, $\hat{F}_n(x) = \int_{-\infty}^{x} \hat{f}_n(r) \, dr$. The first characterization is in terms of $\hat{\varphi}_n$ and perturbation functions.

**Theorem 2.2.** *Let $\widetilde{\varphi}$ be a concave function such that $\{x : \widetilde{\varphi}(x) > -\infty\} = [X_1, X_n]$. Then, $\widetilde{\varphi} = \hat{\varphi}_n$ if and only if*

$$\int \Delta(x) \, d\mathbb{F}_n(x) \leq \int \Delta(x) \exp \widetilde{\varphi}(x) \, dx \tag{2}$$

*for any $\Delta : \mathbb{R} \to \mathbb{R}$ such that $\widetilde{\varphi} + \lambda \Delta$ is concave for some $\lambda > 0$.*

Plugging suitable perturbation functions $\Delta$ in Theorem 2.2 yields valuable information about $\hat{\varphi}_n$ and $\hat{F}_n$. For a first illustration, let $\mu(G)$ and $\mathrm{Var}(G)$ be the mean and variance, respectively, of a distribution (function) $G$ on the real line with finite second moment. Setting $\Delta(x) := \pm x$ or $\Delta(x) := -x^2$ in Theorem 2.4 yields the following.

**Corollary 2.3.**

$$\mu(\hat{F}_n) = \mu(\mathbb{F}_n) \quad \text{and} \quad \mathrm{Var}(\hat{F}_n) \leq \mathrm{Var}(\mathbb{F}_n).$$

Our second characterization is in terms of the empirical distribution function $\mathbb{F}_n$ and the estimated distribution function $\hat{F}_n$. For a continuous and piecewise linear function $h : [X_1, X_n] \to \mathbb{R}$, we define the set of its "knots" to be

$$\mathcal{S}_n(h) := \{t \in (X_1, X_n) : h'(t-) \neq h'(t+)\} \cup \{X_1, X_n\}.$$

Recall that $\hat{\varphi}_n$ is an example of such a function $h$ with $\mathcal{S}_n(\hat{\varphi}_n) \subset \{X_1, X_2, \ldots, X_n\}$.

**Theorem 2.4.** *Let $\widetilde{\varphi}$ be a concave function which is linear on all intervals $[X_j, X_{j+1}]$, $1 \leq j < n$, while $\widetilde{\varphi} = -\infty$ on $\mathbb{R} \setminus [X_1, X_n]$. Defining $\widetilde{F}(x) := \int_{-\infty}^{x} \exp \widetilde{\varphi}(r) \, dr$, we assume further that $\widetilde{F}(X_n) = 1$. Then, $\widetilde{\varphi} = \hat{\varphi}_n$ and $\widetilde{F} = \hat{F}_n$ if, and only if for arbitrary $t \in [X_1, X_n]$,*

$$\int_{X_1}^{t} \widetilde{F}(r) \, dr \leq \int_{X_1}^{t} \mathbb{F}_n(r) \, dr \tag{3}$$

*with equality in the case of $t \in \mathcal{S}_n(\widetilde{\varphi})$.*

A particular consequence of Theorem 2.4 is that the distribution function estimator $\hat{F}_n$ is very close to the empirical distribution function $\mathbb{F}_n$ on $\mathcal{S}_n(\hat{\varphi}_n)$.

**Corollary 2.5.**

$$\mathbb{F}_n - n^{-1} \leq \hat{F}_n \leq \mathbb{F}_n \quad \text{on } \mathcal{S}_n(\hat{\varphi}_n).$$

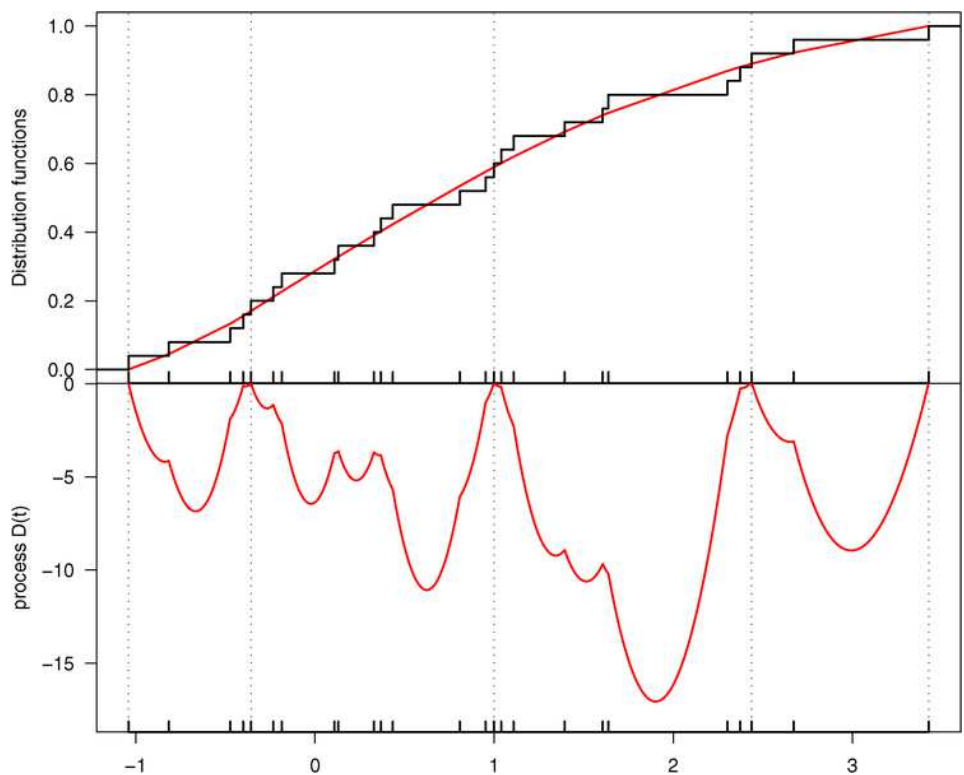

**Figure 1.** Distribution functions and the process $D(t)$ for a Gumbel sample.

Figure 1 illustrates Theorem 2.4 and Corollary 2.5. The upper plot displays $\mathbb{F}_n$ and $\hat{F}_n$ for a sample of $n = 25$ random numbers generated from a Gumbel distribution with density $f(x) = \mathrm{e}^{-x} \exp(-\mathrm{e}^{-x})$ on $\mathbb{R}$. The dotted vertical lines indicate the "kinks" of $\hat{\varphi}_n$, that is, all $t \in \mathcal{S}_n(\hat{\varphi}_n)$. Note that $\hat{F}_n$ and $\mathbb{F}_n$ are indeed very close on the latter set, with equality at the right end-point $X_n$. The lower plot shows the process

$$D(t) := \int_{X_1}^{t} (\hat{F}_n - \mathbb{F}_n)(r) \, \mathrm{d}r$$

for $t \in [X_1, X_n]$. As predicted by Theorem 2.4, this process is non-positive and equals zero on $\mathcal{S}_n(\hat{\varphi}_n)$.

## 3. A data example

In a recent consulting case, a company asked for Monte Carlo experiments to predict the reliability of a certain device they produce. The reliability depended in a certain deterministic way on five different and independent random input parameters. For each input parameter, a sample was available and the goal was to fit a suitable distribution to simulate from. Here, we focus on just one of these input parameters.

At first, we considered two standard approaches to estimate the unknown density $f$, namely, (i) fitting a Gaussian density $\hat{f}_{\mathrm{par}}$ with mean $\mu(\mathbb{F}_n)$ and variance $\hat{\sigma}^2 := n(n-1)^{-1}\mathrm{Var}(\mathbb{F}_n)$; (ii) the kernel density estimator

$$\hat{f}_{\mathrm{ker}}(x) := \int \phi_{\hat{\sigma}/\sqrt{n}}(x-y)\,\mathrm{d}\mathbb{F}_n(y),$$

where $\phi_\sigma$ denotes the density of $\mathcal{N}(0,\sigma^2)$. This very small bandwidth $\hat{\sigma}/\sqrt{n}$ was chosen to obtain a density with variance $\hat{\sigma}^2$ and to avoid putting too much weight into the tails.

Looking at the data, approach (i) is clearly inappropriate because our sample of size $n = 787$ revealed a skewed and significantly non-Gaussian distribution. This can be seen in Figure 2, where the multimodal curve corresponds to $\hat{f}_{\mathrm{ker}}$, while the dashed line depicts $\hat{f}_{\mathrm{par}}$. Approach (ii) yielded Monte Carlo results agreeing well with measured reliabilities, but the engineers questioned the multimodality of $\hat{f}_{\mathrm{ker}}$. Choosing a kernel estimator with larger bandwidth would overestimate the variance and put too much weight into the tails. Thus, we agreed on a third approach and estimated $f$ by a slightly smoothed version of $\hat{f}_n$,

$$\hat{f}_n^* := \int \phi_{\hat{\gamma}}(x-y)\,\mathrm{d}\hat{F}_n(y),$$

with $\hat{\gamma}^2 := \hat{\sigma}^2 - \mathrm{Var}(\hat{F}_n)$, so that the variance of $\hat{f}_n^*$ coincides with $\hat{\sigma}^2$. Since log-concavity is preserved under convolution (cf. Prékopa (1971)), $\hat{f}_n^*$ is also log-concave. For the explicit computation of $\mathrm{Var}(\hat{F}_n)$, see Dümbgen et al. (2007a). By smoothing, we also avoid the small discontinuities of $\hat{f}_n$ at $X_1$ and $X_n$. This density estimator is the skewed unimodal curve in Figure 2. It also yielded convincing results in the Monte Carlo simulations.

Note that both estimators $\hat{f}_n$ and $\hat{f}_n^*$ are fully automatic. Moreover, it is very easy to sample from these densities: let $\mathcal{S}_n(\hat{\varphi}_n)$ consist of $x_0 < x_1 < \cdots < x_m$, and consider the data $X_i$ temporarily as fixed. Now,

(a) generate a random index $J \in \{1,2,\ldots,m\}$ with $\mathbb{P}(J=j) = \hat{F}_n(x_j) - \hat{F}_n(x_{j-1})$;
(b) generate

$$X := x_{J-1} + (x_J - x_{J-1}) \cdot \begin{cases} \log(1+(\mathrm{e}^\Theta-1)U)/\Theta, & \text{if } \Theta \neq 0, \\ U, & \text{if } \Theta = 0, \end{cases}$$

where $\Theta := \hat{\varphi}_n(x_J) - \hat{\varphi}_n(x_{J-1})$ and $U \sim \mathrm{Unif}[0,1]$;
(c) generate

$$X^* := X + \hat{\gamma}Z \quad \text{with } Z \sim \mathcal{N}(0,1),$$

where $J$, $U$ and $Z$ are independent. Then, $X \sim \hat{f}_n$ and $X^* \sim \hat{f}_n^*$.

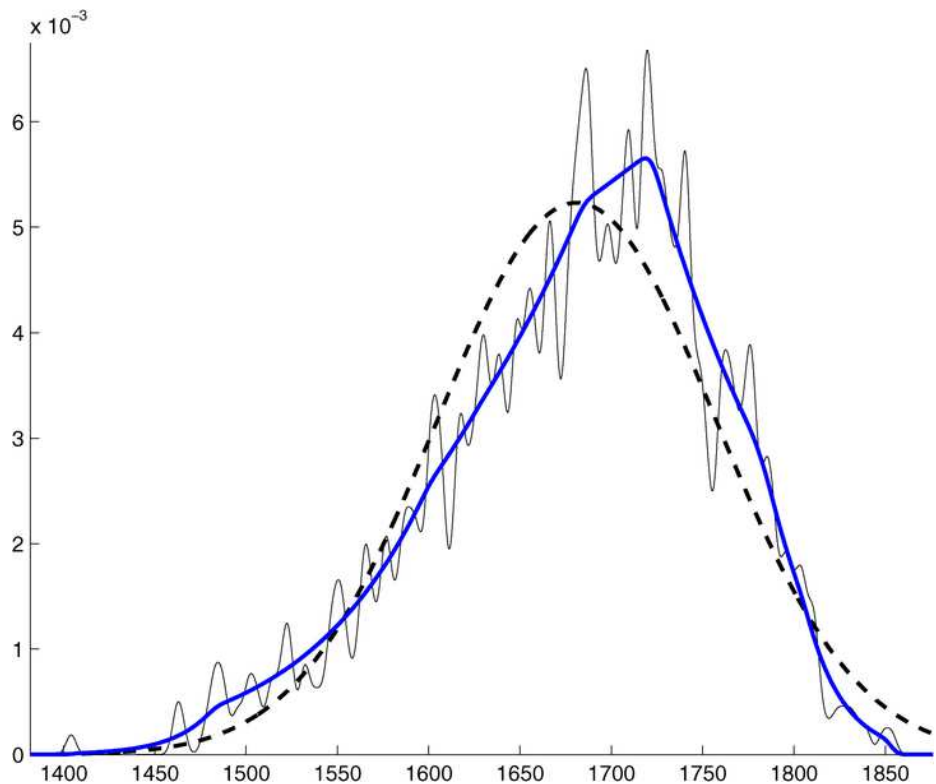

**Figure 2.** Three competing density estimators.

## 4. Uniform consistency

Let us introduce some notation. For any integer $n > 1$, we define

$$\rho_n := \log(n)/n$$

and the uniform norm of a function $g : I \to \mathbb{R}$ on an interval $I \subset \mathbb{R}$ is denoted by

$$\|g\|_\infty^I := \sup_{x \in I} |g(x)|.$$

We say that $g$ belongs to the Hölder class $\mathcal{H}^{\beta,L}(I)$ with exponent $\beta \in [1,2]$ and constant $L > 0$ if for all $x, y \in I$, we have

$$|g(x) - g(y)| \leq L|x - y|, \qquad \text{if } \beta = 1,$$
$$|g'(x) - g'(y)| \leq L|x - y|^{\beta - 1}, \qquad \text{if } \beta > 1.$$

*Uniform consistency of $\hat{\varphi}_n$.* Our main result is the following theorem.

**Theorem 4.1.** *Assume for the log-density $\varphi = \log f$ that $\varphi \in \mathcal{H}^{\beta,L}(T)$ for some exponent $\beta \in [1,2]$, some constant $L > 0$ and a subinterval $T = [A,B]$ of the interior of $\{f > 0\}$. Then,*

$$\max_{t \in T}(\hat{\varphi}_n - \varphi)(t) = O_{\mathrm{p}}(\rho_n^{\beta/(2\beta+1)}),$$

$$\max_{t \in T(n,\beta)}(\varphi - \hat{\varphi}_n)(t) = O_{\mathrm{p}}(\rho_n^{\beta/(2\beta+1)}),$$

*where $T(n,\beta) := [A + \rho_n^{1/(2\beta+1)}, B - \rho_n^{1/(2\beta+1)}]$.*

Note that the previous result remains true when we replace $\hat{\varphi}_n - \varphi$ with $\hat{f}_n - f$. It is well known that the rates of convergence in Theorem 4.1 are optimal, even if $\beta$ was known (cf. Khas'minskii (1978)). Thus, our estimators adapt to the unknown smoothness of $f$ in the range $\beta \in [1,2]$.

Also, note that concavity of $\varphi$ implies that it is Lipschitz-continuous, that is, belongs to $\mathcal{H}^{1,L}(T)$ for some $L > 0$ on any interval $T = [A,B]$ with $A > \inf\{f > 0\}$ and $B < \sup\{f > 0\}$. Hence, one can easily deduce from Theorem 4.1 that $\hat{f}_n$ is consistent in $L_1(\mathbb{R})$ and that $\hat{F}_n$ is uniformly consistent.

**Corollary 4.2.**

$$\int |\hat{f}_n(x) - f(x)|\,\mathrm{d}x \to_{\mathrm{p}} 0 \quad \textit{and} \quad \|\hat{F}_n - F\|_\infty^{\mathbb{R}} \to_{\mathrm{p}} 0.$$

*Distance of two consecutive knots and uniform consistency of $\hat{F}_n$.* By means of Theorem 4.1, we can solve a "gap problem" for log-concave density estimation. The term "gap problem" was first used by Balabdaoui and Wellner (2008) to describe the problem of computing the distance between two consecutive knots of certain estimators.

**Theorem 4.3.** *Suppose that the assumptions of Theorem 4.1 hold. Assume, further, that $\varphi'(x) - \varphi'(y) \geq C(y - x)$ for some constant $C > 0$ and arbitrary $A \leq x < y \leq B$, where $\varphi'$ stands for $\varphi'(\cdot-)$ or $\varphi'(\cdot+)$. Then,*

$$\sup_{x \in T} \min_{y \in \mathcal{S}_n(\hat{\varphi}_n)} |x - y| = O_{\mathrm{p}}(\rho_n^{\beta/(4\beta+2)}).$$

Theorems 4.1 and 4.3, combined with a result of Stute (1982) about the modulus of continuity of empirical processes, yield a rate of convergence for the maximal difference between $\hat{F}_n$ and $\mathbb{F}_n$ on compact intervals.

**Theorem 4.4.** *Under the assumptions of Theorem 4.3,*

$$\max_{t \in T(n,\beta)}|\hat{F}_n(t) - \mathbb{F}_n(t)| = O_{\mathrm{p}}(\rho_n^{3\beta/(4\beta+2)}).$$

*In particular, if $\beta > 1$, then*

$$\max_{t \in T(n,\beta)} |\hat{F}_n(t) - \mathbb{F}_n(t)| = o_{\mathrm{p}}(n^{-1/2}).$$

Thus, under certain regularity conditions, the estimators $\hat{F}_n$ and $\mathbb{F}_n$ are asymptotically equivalent on compact sets. Conclusions of this type are known for the Grenander estimator (cf. Kiefer and Wolfowitz (1976)) and the least squares estimator of a convex density on $[0, \infty)$ (cf. Balabdaoui and Wellner (2007)).

The result of Theorem 4.4 is also related to recent results of Giné and Nickl (2007, 2008). In the latter paper, they devise kernel density estimators with data-driven bandwidths which are also adaptive with respect to $\beta$ in a certain range, while the integrated density estimator is asymptotically equivalent to $\mathbb{F}_n$ on the whole real line. However, if $\beta \geq 3/2$, they must use kernel functions of higher order, that is, no longer non-negative, and simulating data from the resulting estimated density is not straightforward.

*Example.* Let us illustrate Theorems 4.1 and 4.4 with simulated data, again from the Gumbel distribution with $\varphi(x) = -x - \mathrm{e}^{-x}$. Here, $\varphi''(x) = -\mathrm{e}^{-x}$, so the assumptions of our theorems are satisfied with $\beta = 2$ for any compact interval $T$. The upper panels of Figure 3 show the true log-density $\varphi$ (dashed line) and the estimator $\hat{\varphi}_n$ (line) for samples of sizes $n = 200$ (left) and $n = 2000$ (right). The lower panels show the corresponding empirical processes $n^{1/2}(\mathbb{F}_n - F)$ (jagged curves) and $n^{1/2}(\hat{F}_n - F)$ (smooth curves). First, the quality of the estimator $\hat{\varphi}_n$ is quite good, even in the tails, and the quality increases with sample size, as expected. Looking at the empirical processes, the similarity between $n^{1/2}(\mathbb{F}_n - F)$ and $n^{1/2}(\hat{F}_n - F)$ increases with sample size, too, but rather slowly. Also, note that the estimator $\hat{F}_n$ outperforms $\mathbb{F}_n$ in terms of supremum distance from $F$, which leads us to the next paragraph.

*Marshall's lemma.* In all simulations we looked at, the estimator $\hat{F}_n$ satisfied the inequality

$$\|\hat{F}_n - F\|_\infty^{\mathbb{R}} \leq \|\mathbb{F}_n - F\|_\infty^{\mathbb{R}}, \tag{4}$$

provided that $f$ is indeed log-concave. Figure 3 shows two numerical examples of this phenomenon. In view of such examples and Marshall's (1970) lemma about the Grenander estimator $\hat{F}_n^{\mathrm{mon}}$, we first tried to verify that (4) is correct almost surely and for any $n > 1$. However, one can construct counterexamples showing that (4) may be violated, even if the right-hand side is multiplied with any fixed constant $C > 1$. Nevertheless, our first attempts resulted in a version of Marshall's lemma for *convex* density estimation; see Dümbgen *et al.* (2007). For the present setting, we conjecture that (4) is true with asymptotic probability one as $n \to \infty$, that is,

$$\mathbb{P}(\|\hat{F}_n - F\|_\infty^{\mathbb{R}} \leq \|\mathbb{F}_n - F\|_\infty^{\mathbb{R}}) \to 1.$$

*A monotone hazard rate estimator.* Estimation of a monotone hazard rate is described, for instance, in the book by Robertson *et al.* (1988). They directly solve an isotonic estimation problem similar to that for the Grenander density estimator. For this setting, Hall *et al.* (2001) and Hall and van Keilegom (2005) consider methods based upon suitable modifications of kernel estimators. Alternatively, in our setting, it follows from Lemma A.2 in Section 5 that

$$\hat{h}_n(x) := \frac{\hat{f}_n(x)}{1 - \hat{F}_n(x)}$$

defines a simple plug-in estimator of the hazard rate on $(-\infty, X_n)$ which is also nondecreasing. By virtue of Theorem 4.1 and Corollary 4.2, it is uniformly consistent on any compact subinterval of the interior of $\{f > 0\}$. Theorems 4.1 and 4.4 even entail a rate of convergence, as follows.

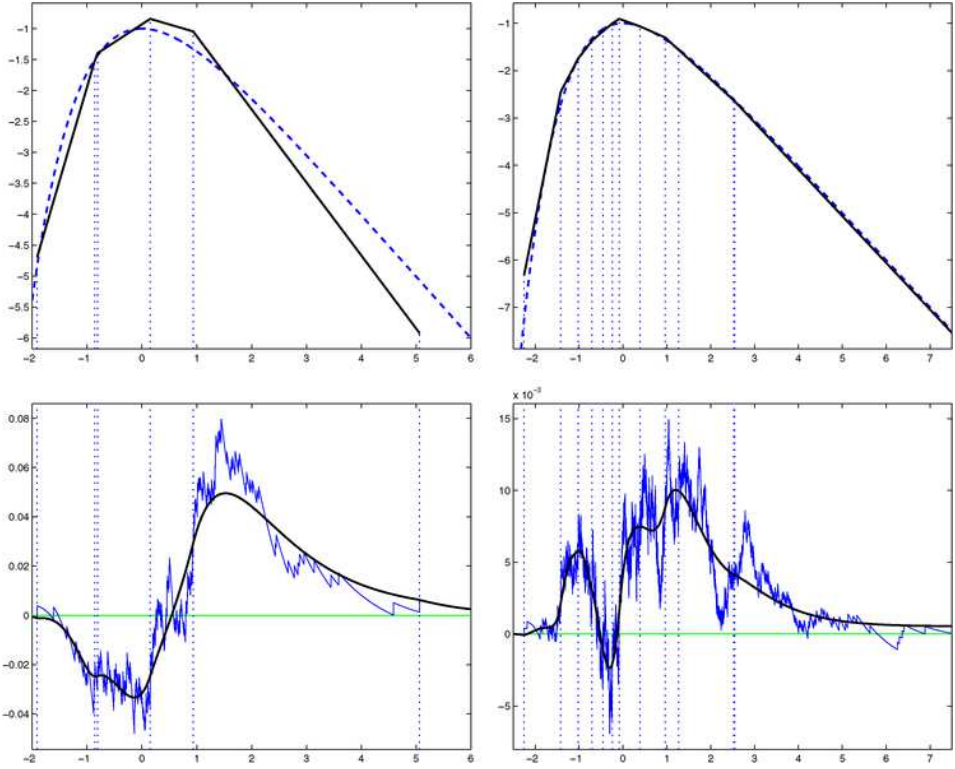

**Figure 3.** Density functions and empirical processes for Gumbel samples of size $n = 200$ and $n = 2000$.

**Corollary 4.5.** *Under the assumptions of Theorem 4.3,*

$$\max_{t \in T(n,\beta)} |\hat{h}_n(t) - h(t)| = O_p(\rho_n^{\beta/(2\beta+1)}).$$

## 5. Outlook

Starting from the results presented here, Balabdaoui *et al.* (2008) recently derived the pointwise limiting distribution of $\hat{f}_n$. They also considered the limiting distribution of $\mathrm{argmax}_{x \in \mathbb{R}} \hat{f}_n(x)$ as an estimator of the mode of $f$. Empirical findings of Müller and Rufibach (2008) show that the estimator $\hat{f}_n$ is even useful for extreme value statistics. Log-concave densities also have potential as building blocks in more complex models (e.g., regression or classification) or when handling censored data (cf. Dümbgen *et al.* (2007a)).

Unfortunately, our proofs work only for fixed compact intervals, whereas simulations suggest that the estimators perform well on the whole real line. Presently, the authors are working on a different approach, where $\hat{\varphi}_n$ is represented locally as a parametric maximum likelihood estimator of a log-linear density. Presumably, this will deepen our understanding of the log-concave NPMLE's consistency properties, particularly in the tails. For instance, we conjecture that $\mathbb{F}_n$ and $\hat{F}_n$ are asymptotically equivalent on any interval $T$ on which $\varphi'$ is strictly decreasing.

## Appendix: Auxiliary results and proofs

### A.1. Two facts about log-concave densities

The following two results about a log-concave density $f = \exp \varphi$ and its distribution function $F$ are of independent interest. The first result entails that the density $f$ has at least subexponential tails.

**Lemma A.1.** *For arbitrary points* $x_1 < x_2$,

$$\sqrt{f(x_1)f(x_2)} \leq \frac{F(x_2) - F(x_1)}{x_2 - x_1}.$$

*Moreover, for* $x_o \in \{f > 0\}$ *and any real* $x \neq x_o$,

$$\frac{f(x)}{f(x_o)} \leq \begin{cases} \left( \dfrac{h(x_o, x)}{f(x_o)|x - x_o|} \right)^2, \\ \exp\left( 1 - \dfrac{f(x_o)|x - x_o|}{h(x_o, x)} \right) & \text{if } f(x_o)|x - x_o| \geq h(x_o, x), \end{cases}$$

*where*

$$h(x_o, x) := F(\max(x_o, x)) - F(\min(x_o, x)) \leq \begin{cases} F(x_o), & \text{if } x < x_o, \\ 1 - F(x_o), & \text{if } x > x_o. \end{cases}$$

A second well-known result (Barlow and Proschan (1975), Lemma 5.8) provides further connections between the density $f$ and the distribution function $F$. In particular, it entails that $f/(F(1-F))$ is bounded away from zero on $\{x : 0 < F(x) < 1\}$.

**Lemma A.2.** *The function $f/F$ is non-increasing on $\{x : 0 < F(x) \leq 1\}$ and the function $f/(1-F)$ is non-decreasing on $\{x : 0 \leq F(x) < 1\}$.*

**Proof of Lemma A.1.** To prove the first inequality, it suffices to consider the non-trivial case of $x_1, x_2 \in \{f > 0\}$. Concavity of $\varphi$ then entails that

$$F(x_2) - F(x_1) \geq \int_{x_1}^{x_2} \exp\left(\frac{x_2 - t}{x_2 - x_1}\varphi(x_1) + \frac{t - x_1}{x_2 - x_1}\varphi(x_2)\right) dt$$

$$= (x_2 - x_1) \int_0^1 \exp((1-u)\varphi(x_1) + u\varphi(x_2)) \, du$$

$$\geq (x_2 - x_1) \exp\left(\int_0^1 ((1-u)\varphi(x_1) + u\varphi(x_2)) \, du\right)$$

$$= (x_2 - x_1) \exp(\varphi(x_1)/2 + \varphi(x_2)/2)$$

$$= (x_2 - x_1)\sqrt{f(x_1)f(x_2)},$$

where the second inequality follows from Jensen's inequality.

We prove the second asserted inequality only for $x > x_o$, that is, $h(x_o, x) = F(x) - F(x_o)$, the other case being handled analogously. The first part entails that

$$\frac{f(x)}{f(x_o)} \leq \left(\frac{h(x_o, x)}{f(x_o)(x - x_o)}\right)^2,$$

and the right-hand side is not greater than one if $f(x_o)(x - x_o) \geq h(x_o, x)$. In the latter case, recall that

$$h(x_o, x) \geq (x - x_o) \int_0^1 \exp((1-u)\varphi(x_o) + u\varphi(x)) \, du = f(x_o)(x - x_o)J(\varphi(x) - \varphi(x_o))$$

with $\varphi(x) - \varphi(x_o) \leq 0$, where $J(y) := \int_0^1 \exp(uy) \, du$. Elementary calculations show that $J(-r) = (1 - e^{-r})/r \geq 1/(1 + r)$ for arbitrary $r > 0$. Thus,

$$h(x_o, x) \geq \frac{f(x_o)(x - x_o)}{1 + \varphi(x_o) - \varphi(x)},$$

which is equivalent to $f(x)/f(x_o) \leq \exp(1 - f(x_o)(x - x_o)/h(x_o, x))$. □

### A.2. Proofs of the characterizations

**Proof of Theorem 2.2.** In view of Theorem 2.1, we may restrict our attention to concave and real-valued functions $\varphi$ on $[X_1, X_n]$ and set $\varphi := -\infty$ on $\mathbb{R} \setminus [X_1, X_n]$. The set $\mathcal{C}_n$ of all such functions is a convex cone and for any function $\Delta : \mathbb{R} \to \mathbb{R}$ and $t > 0$, concavity of $\varphi + t\Delta$ on $\mathbb{R}$ is equivalent to its concavity on $[X_1, X_n]$.

One can easily verify that $\Psi_n$ is a concave and real-valued functional on $\mathcal{C}_n$. Hence, as well known from convex analysis, a function $\widetilde{\varphi} \in \mathcal{C}_n$ maximizes $\Psi_n$ if and only if

$$\lim_{t \downarrow 0} \frac{\Psi_n(\widetilde{\varphi} + t(\varphi - \widetilde{\varphi})) - \Psi_n(\widetilde{\varphi})}{t} \leq 0$$

for all $\varphi \in \mathcal{C}_n$. But, this is equivalent to the requirement that

$$\lim_{t \downarrow 0} \frac{\Psi_n(\widetilde{\varphi} + t\Delta) - \Psi_n(\widetilde{\varphi})}{t} \leq 0$$

for any function $\Delta : \mathbb{R} \to \mathbb{R}$ such that $\widetilde{\varphi} + \lambda\Delta$ is concave for some $\lambda > 0$. The assertion of the theorem now follows from

$$\lim_{t \downarrow 0} \frac{\Psi_n(\widetilde{\varphi} + t\Delta) - \Psi_n(\widetilde{\varphi})}{t} = \int \Delta \, d\mathbb{F}_n - \int \Delta(x) \exp\widetilde{\varphi}(x) \, dx. \qquad \square$$

**Proof of Theorem 2.4.** We start with a general observation. Let $G$ be some distribution (function) with support $[X_1, X_n]$ and let $\Delta : [X_1, X_n] \to \mathbb{R}$ be absolutely continuous with $L_1$-derivative $\Delta'$. It then follows from Fubini's theorem that

$$\int \Delta \, dG = \Delta(X_n) - \int_{X_1}^{X_n} \Delta'(r) G(r) \, dr. \tag{A.1}$$

Now, suppose that $\widetilde{\varphi} = \hat{\varphi}_n$ and let $t \in (X_1, X_n]$. Let $\Delta$ be absolutely continuous on $[X_1, X_n]$ with $L_1$-derivative $\Delta'(r) = 1\{r \leq t\}$ and arbitrary value of $\Delta(X_n)$. Clearly, $\widetilde{\varphi} + \Delta$ is concave, whence (2) and (A.1) entail that

$$\Delta(X_n) - \int_{X_1}^{t} \mathbb{F}_n(r) \, dr \leq \Delta(X_n) - \int_{X_1}^{t} \widetilde{F}(r) \, dr,$$

which is equivalent to inequality (3). In the case of $t \in \mathcal{S}_n(\widetilde{\varphi}) \setminus \{X_1\}$, let $\Delta'(r) = -1\{r \leq t\}$. Then, $\widetilde{\varphi} + \lambda\Delta$ is concave for some $\lambda > 0$ so that

$$\Delta(X_n) + \int_{X_1}^{t} \mathbb{F}_n(r) \, dr \leq \Delta(X_n) + \int_{X_1}^{t} \widetilde{F}(r) \, dr,$$

which yields equality in (3).

Now, suppose that $\widetilde{\varphi}$ satisfies inequality (3) for all $t$ with equality if $t \in \mathcal{S}_n(\widetilde{\varphi})$. In view of Theorem 2.1 and the proof of Theorem 2.2, it suffices to show that (2) holds for any function $\Delta$ defined on $[X_1, X_n]$ which is linear on each interval $[X_j, X_{j+1}]$, $1 \leq j < n$, while $\widetilde{\varphi} + \lambda\Delta$ is concave for some $\lambda > 0$. The latter requirement is equivalent to $\Delta$ being concave between two consecutive knots of $\widetilde{\varphi}$. Elementary considerations show that the $L_1$-derivative of such a function $\Delta$ may be written as

$$\Delta'(r) = \sum_{j=2}^{n} \beta_j 1\{r \leq X_j\},$$

with real numbers $\beta_2, \ldots, \beta_n$ such that

$$\beta_j \geq 0 \qquad \text{if } X_j \notin \mathcal{S}_n(\widetilde{\varphi}).$$

Consequently, it follows from (A.1) and our assumptions on $\widetilde{\varphi}$ that

$$\int \Delta \, d\mathbb{F}_n = \Delta(X_n) - \sum_{j=2}^{n} \beta_j \int_{X_1}^{X_j} \mathbb{F}_n(r) \, dr$$

$$\leq \Delta(X_n) - \sum_{j=2}^{n} \beta_j \int_{X_1}^{X_j} \widetilde{F}(r) \, dr$$

$$= \int \Delta \, d\widetilde{F}. \qquad \square$$

**Proof of Corollary 2.5.** For $t \in \mathcal{S}_n(\hat{\varphi}_n)$ and $s < t < u$, it follows from Theorem 2.4 that

$$\frac{1}{u-t} \int_t^u \hat{F}_n(r) \, dr \leq \frac{1}{u-t} \int_s^t \mathbb{F}_n(r) \, dr$$

and

$$\frac{1}{t-s} \int_s^t \hat{F}_n(r) \, dr \geq \frac{1}{t-s} \int_s^t \mathbb{F}_n(r) \, dr.$$

Letting $u \downarrow t$ and $s \uparrow t$ yields

$$\hat{F}_n(t) \leq \mathbb{F}_n(t) \quad \text{and} \quad \hat{F}_n(t) \geq \mathbb{F}_n(t-) = \mathbb{F}_n(t) - n^{-1}. \qquad \square$$

## A.3. Proof of $\hat{\varphi}_n$'s consistency

Our proof of Theorem 4.1 involves a refinement and modification of methods introduced by Dümbgen *et al.* (2004). A first key ingredient is an inequality for concave functions due to Dümbgen (1998) (see also Dümbgen *et al.* (2004) or Rufibach (2006)).

**Lemma A.3.** *For any $\beta \in [1,2]$ and $L > 0$, there exists a constant $K = K(\beta, L) \in (0,1]$ with the following property. Suppose that $g$ and $\hat{g}$ are concave and real-valued functions on a compact interval $T = [A, B]$, where $g \in \mathcal{H}^{\beta, L}(T)$. Let $\epsilon > 0$ and $0 < \delta \leq K \min\{B - A, \epsilon^{1/\beta}\}$. Then*

$$\sup_{t \in T} (\hat{g} - g) \geq \epsilon \quad or \quad \sup_{t \in [A+\delta, B-\delta]} (g - \hat{g}) \geq \epsilon$$

*implies that*

$$\inf_{t \in [c, c+\delta]} (\hat{g} - g)(t) \geq \epsilon/4 \quad or \quad \inf_{t \in [c, c+\delta]} (g - \hat{g})(t) \geq \epsilon/4$$

*for some $c \in [A, B - \delta]$.*

Starting from this lemma, let us first sketch the idea of our proof of Theorem 4.1. Suppose we had a family $\mathcal{D}$ of measurable functions $\Delta$ with finite seminorm

$$\sigma(\Delta) := \left(\int \Delta^2 \,\mathrm{d}F\right)^{1/2},$$

such that

$$\sup_{\Delta \in \mathcal{D}} \frac{|\int \Delta \,\mathrm{d}(\mathbb{F}_n - F)|}{\sigma(\Delta)\rho_n^{1/2}} \leq C \tag{A.2}$$

with asymptotic probability one, where $C > 0$ is some constant. If, in addition, $\varphi - \hat{\varphi}_n \in \mathcal{D}$ and $\varphi - \hat{\varphi}_n \leq C$ with asymptotic probability one, then we could conclude that

$$\left|\int (\varphi - \hat{\varphi}_n) \,\mathrm{d}(\mathbb{F}_n - F)\right| \leq C\sigma(\varphi - \hat{\varphi}_n)\rho_n^{1/2},$$

while Theorem 2.2, applied to $\Delta := \varphi - \hat{\varphi}_n$, entails that

$$\int (\varphi - \hat{\varphi}_n) \,\mathrm{d}(\mathbb{F}_n - F) \leq \int (\varphi - \hat{\varphi}_n) \,\mathrm{d}(\hat{F} - F)$$

$$= -\int \Delta(1 - \exp(-\Delta)) \,\mathrm{d}F$$

$$\leq -(1+C)^{-1} \int \Delta^2 \,\mathrm{d}F$$

$$= -(1+C)^{-1}\sigma(\varphi - \hat{\varphi}_n)^2$$

because $y(1 - \exp(-y)) \geq (1 + y_+)^{-1} y^2$ for all real $y$, where $y_+ := \max(y, 0)$. Hence, with asymptotic probability one,

$$\sigma(\varphi - \hat{\varphi}_n)^2 \leq C^2(1+C)^2 \rho_n.$$

Now, suppose that $|\varphi - \hat{\varphi}_n| \geq \epsilon_n$ on a subinterval of $T = [A, B]$ of length $\epsilon_n^{1/\beta}$, where $(\epsilon_n)_n$ is a fixed sequence of numbers $\epsilon_n > 0$ tending to zero. Then, $\sigma(\varphi - \hat{\varphi}_n)^2 \geq \epsilon_n^{(2\beta+1)/\beta} \min_T(f)$, so that

$$\epsilon_n \leq \widetilde{C} \rho_n^{2\beta/(2\beta+1)}$$

with $\widetilde{C} = (C^2(1+C)^2/\min_T(f))^{\beta/(2\beta+1)}$.

The previous considerations will be modified in two aspects to get a rigorous proof of Theorem 4.1. For technical reasons, we must replace the denominator $\sigma(\Delta)\rho_n^{1/2}$ of inequality (A.2) with $\sigma(\Delta)\rho_n^{1/2} + W(\Delta)\rho_n^{2/3}$, where

$$W(\Delta) := \sup_{x \in \mathbb{R}} \frac{|\Delta(x)|}{\max(1, |\varphi(x)|)}.$$

This is necessary to deal with functions $\Delta$ with small values of $F(\{\Delta \neq 0\})$. Moreover, we shall work with simple "caricatures" of $\varphi - \hat{\varphi}_n$, namely, functions which are piecewise linear with at most three knots. Throughout this section, piecewise linearity does *not* necessarily imply continuity. A function being piecewise linear with at most $m$ knots means that the real line may be partitioned into $m+1$ non-degenerate intervals on each of which the function is linear. Then, the $m$ real boundary points of these intervals are the knots.

The next lemma extends inequality (2) to certain piecewise linear functions.

**Lemma A.4.** *Let $\Delta : \mathbb{R} \to \mathbb{R}$ be piecewise linear such that each knot $q$ of $\Delta$ satisfies one of the following two properties:*

$$q \in \mathcal{S}_n(\hat{\varphi}_n) \quad and \quad \Delta(q) = \liminf_{x \to q} \Delta(x); \tag{A.3}$$

$$\Delta(q) = \lim_{r \to q} \Delta(r) \quad and \quad \Delta'(q-) \geq \Delta'(q+). \tag{A.4}$$

*Then,*

$$\int \Delta \, \mathrm{d}\mathbb{F}_n \leq \int \Delta \, \mathrm{d}\hat{F}_n. \tag{A.5}$$

We can now specify the "caricatures" mentioned above.

**Lemma A.5.** *Let $T = [A, B]$ be a fixed subinterval of the interior of $\{f > 0\}$. Let $\varphi - \hat{\varphi}_n \geq \epsilon$ or $\hat{\varphi}_n - \varphi \geq \epsilon$ on some interval $[c, c + \delta] \subset T$ with length $\delta > 0$ and suppose that $X_1 < c$ and $X_n > c + \delta$. There then exists a piecewise linear function $\Delta$ with at most three knots, each of which satisfies condition (A.3) or (A.4), and a positive constant $K' = K'(f, T)$ such that*

$$|\varphi - \hat{\varphi}_n| \geq \epsilon |\Delta|, \tag{A.6}$$

$$\Delta(\varphi - \hat{\varphi}_n) \geq 0, \tag{A.7}$$

$$\Delta \leq 1, \tag{A.8}$$

$$\int_c^{c+\delta} \Delta^2(x) \, \mathrm{d}x \geq \delta/3, \tag{A.9}$$

$$W(\Delta) \leq K' \delta^{-1/2} \sigma(\Delta). \tag{A.10}$$

Our last ingredient is a surrogate for (A.2).

**Lemma A.6.** *Let $\mathcal{D}_m$ be the family of all piecewise linear functions on $\mathbb{R}$ with at most $m$ knots. There exists a constant $K'' = K''(f)$ such that*

$$\sup_{m \geq 1, \Delta \in \mathcal{D}_m} \frac{|\int \Delta \, \mathrm{d}(\mathbb{F}_n - F)|}{\sigma(\Delta) m^{1/2} \rho_n^{1/2} + W(\Delta) m \rho_n^{2/3}} \leq K'',$$

*with probability tending to one as $n \to \infty$.*

Before we verify all of these auxiliary results, let us proceed with the main proof.

**Proof of Theorem 4.1.** Suppose that

$$\sup_{t \in T}(\hat{\varphi}_n - \varphi)(t) \geq C\epsilon_n$$

or

$$\sup_{t \in [A+\delta_n, B-\delta_n]}(\varphi - \hat{\varphi}_n)(t) \geq C\epsilon_n$$

for some constant $C > 0$, where $\epsilon_n := \rho_n^{\beta/(2\beta+1)}$ and $\delta_n := \rho_n^{1/(2\beta+1)} = \epsilon_n^{1/\beta}$. It follows from Lemma A.3 with $\epsilon := C\epsilon_n$ that in the case of $C \geq K^{-\beta}$ and for sufficiently large $n$, there is a (random) interval $[c_n, c_n + \delta_n] \subset T$ on which either $\hat{\varphi}_n - \varphi \geq (C/4)\epsilon_n$ or $\varphi - \hat{\varphi}_n \geq (C/4)\epsilon_n$. But, then, there is a (random) function $\Delta_n \in \mathcal{D}_3$ fulfilling the conditions stated in Lemma A.5. For this $\Delta_n$, it follows from (A.5) that

$$\int_{\mathbb{R}} \Delta_n \, \mathrm{d}(F - \mathbb{F}_n) \geq \int_{\mathbb{R}} \Delta_n \, \mathrm{d}(F - \hat{F}_n) = \int_{\mathbb{R}} \Delta_n (1 - \exp[-(\varphi - \hat{\varphi}_n)]) \, \mathrm{d}F. \qquad (\text{A.11})$$

With $\widetilde{\Delta}_n := (C/4)\epsilon_n \Delta_n$, it follows from (A.6–A.7) that the right-hand side of (A.11) is not smaller than

$$(4/C)\epsilon_n^{-1} \int \widetilde{\Delta}_n (1 - \exp(-\widetilde{\Delta}_n)) \, \mathrm{d}F \geq \frac{(4/C)\epsilon_n^{-1}}{1 + (C/4)\epsilon_n}\sigma(\widetilde{\Delta}_n)^2 = \frac{(C/4)\epsilon_n}{1 + o(1)}\sigma(\Delta_n)^2$$

because $\widetilde{\Delta}_n \leq (C/4)\epsilon_n$, by (A.8). On the other hand, according to Lemma A.6, we may assume that

$$\int_{\mathbb{R}} \Delta_n \, \mathrm{d}(F - \mathbb{F}_n) \leq K''(3^{1/2}\sigma(\Delta_n)\rho_n^{1/2} + 3W(\Delta_n)\rho_n^{2/3})$$

$$\leq K''(3^{1/2}\rho_n^{1/2} + 3K'\delta_n^{-1/2}\rho_n^{2/3})\sigma(\Delta_n) \qquad (\text{by (A.10)})$$

$$\leq K''(3^{1/2}\rho_n^{1/2} + 3K'\rho_n^{2/3-1/(4\beta+2)})\sigma(\Delta_n)$$

$$\leq G\rho_n^{1/2}\sigma(\Delta_n)$$

for some constant $G = G(\beta, L, f, T)$ because $2/3 - 1/(4\beta+2) \geq 2/3 - 1/6 = 1/2$. Consequently,

$$C^2 \leq \frac{16G^2(1+o(1))\epsilon_n^{-2}\rho_n}{\sigma(\Delta_n)^2} = \frac{16G^2(1+o(1))}{\delta_n^{-1}\sigma(\Delta_n)^2} \leq \frac{48G^2(1+o(1))}{\min_T(f)},$$

where the last inequality follows from (A.9). $\qquad\square$

**Proof of Lemma A.4.** There is a sequence of continuous, piecewise linear functions $\Delta_k$ converging pointwise isotonically to $\Delta$ as $k \to \infty$ such that any knot $q$ of $\Delta_k$ either belongs to $\mathcal{S}_n(\hat{\varphi}_n)$ or $\Delta'_k(q-) > \Delta'_k(q+)$. Thus, $\hat{\varphi}_n + \lambda \Delta_k$ is concave for sufficiently small $\lambda > 0$. Consequently, since $\Delta_1 \leq \Delta_k \leq \Delta$ for all $k$, it follows from dominated convergence and (2) that

$$\int \Delta \, d\mathbb{F}_n = \lim_{k \to \infty} \int \Delta_k \, d\mathbb{F}_n \leq \lim_{k \to \infty} \int \Delta_k \, d\hat{F}_n = \int \Delta \, d\hat{F}_n. \qquad \square$$

**Proof of Lemma A.5.** The crucial point in all the cases we must distinguish is to construct a $\Delta \in \mathcal{D}_3$ satisfying the assumptions of Lemma A.4 and (A.6–A.9). Recall that $\hat{\varphi}_n$ is piecewise linear.

**Case 1a**: $\hat{\varphi}_n - \varphi \geq \epsilon$ on $[c, c+\delta]$ and $\mathcal{S}_n(\hat{\varphi}_n) \cap (c, c+\delta) \neq \varnothing$. Here, we choose a continuous function $\Delta \in \mathcal{D}_3$ with knots $c$, $c+\delta$ and $x_o \in \mathcal{S}_n(\hat{\varphi}_n) \cap (c, c+\delta)$, where $\Delta := 0$ on $(-\infty, c] \cup [c+\delta, \infty)$ and $\Delta(x_o) := -1$. Here, the assumptions of Lemma A.4 and requirements (A.6–A.9) are easily verified.

**Case 1b**: $\hat{\varphi}_n - \varphi \geq \epsilon$ on $[c, c+\delta]$ and $\mathcal{S}_n(\hat{\varphi}_n) \cap (c, c+\delta) = \varnothing$. Let $[c_o, d_o] \supset [c, c+\delta]$ be the maximal interval on which $\varphi - \hat{\varphi}_n$ is concave. There then exists a linear function $\widetilde{\Delta}$ such that $\widetilde{\Delta} \geq \varphi - \hat{\varphi}_n$ on $[c_o, d_o]$ and $\widetilde{\Delta} \leq -\epsilon$ on $[c, c+\delta]$. Next, let $(c_1, d_1) := \{\widetilde{\Delta} < 0\} \cap (c_o, d_o)$. We now define $\Delta \in \mathcal{D}_2$ via

$$\Delta(x) := \begin{cases} 0, & \text{if } x \in (-\infty, c_1) \cup (d_1, \infty), \\ \widetilde{\Delta}/\epsilon, & \text{if } x \in [c_1, d_1]. \end{cases}$$

Again, the assumptions of Lemma A.4 and requirements (A.6–A.9) are easily verified; this time, we even know that $\Delta \leq -1$ on $[c, c+\delta]$, whence $\int_c^{c+\delta} \Delta(x)^2 \, dx \geq \delta$. Figure 4 illustrates this construction.

**Case 2**: $\varphi - \hat{\varphi}_n \geq \epsilon$ on $[c, c+\delta]$. Let $[c_o, c]$ and $[c+\delta, d_o]$ be maximal intervals on which $\hat{\varphi}_n$ is linear. We then define

$$\Delta(x) := \begin{cases} 0, & \text{if } x \in (-\infty, c_o) \cup (d_o, \infty), \\ 1 + \beta_1(x - x_o), & \text{if } x \in [c_o, x_o], \\ 1 + \beta_2(x - x_o), & \text{if } x \in [x_o, d_o], \end{cases}$$

where $x_o := c + \delta/2$ and $\beta_1 \geq 0$ is chosen such that either

$$\Delta(c_o) = 0 \quad \text{and} \quad (\varphi - \hat{\varphi}_n)(c_o) \geq 0, \quad \text{or}$$

$$(\varphi - \hat{\varphi}_n)(c_o) < 0 \quad \text{and} \quad \text{sign}(\Delta) = \text{sign}(\varphi - \hat{\varphi}_n) \qquad \text{on } [c_o, x_o].$$

Analogously, $\beta_2 \leq 0$ is chosen such that

$$\Delta(d_o) = 0 \quad \text{and} \quad (\varphi - \hat{\varphi}_n)(d_o) \geq 0, \quad \text{or}$$

$$(\varphi - \hat{\varphi}_n)(d_o) < 0 \quad \text{and} \quad \text{sign}(\Delta) = \text{sign}(\varphi - \hat{\varphi}_n) \qquad \text{on } [x_o, d_o].$$

Again, the assumptions of Lemma A.4 and requirements (A.6–A.9) are easily verified. Figure 5 depicts an example.

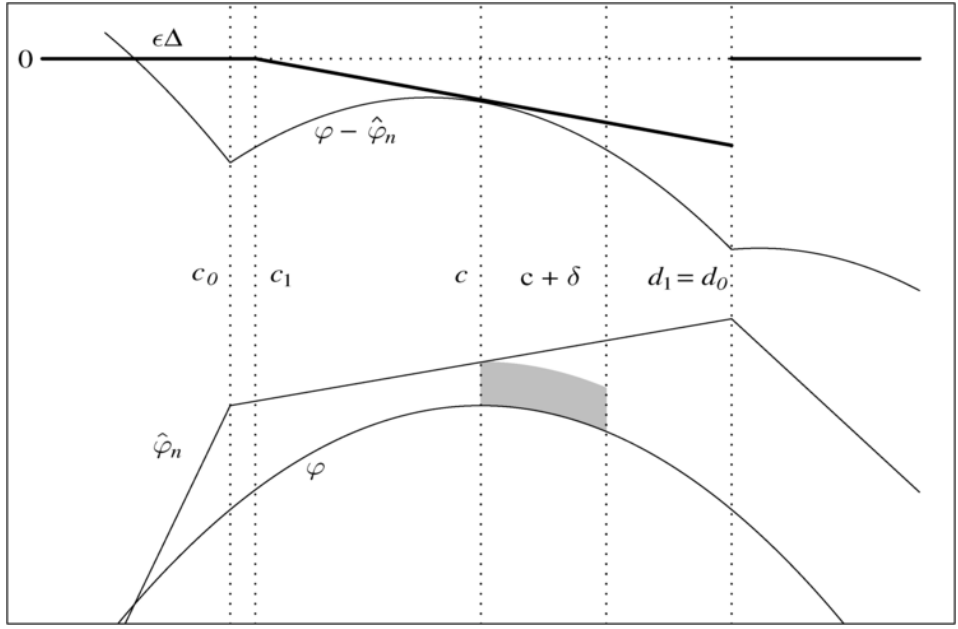

**Figure 4.** The perturbation function $\Delta$ in Case 1b.

It remains to verify requirement (A.10) for our particular functions $\Delta$. Note that by our assumption on $T = [A,B]$, there exist numbers $\tau, C_o > 0$ such that $f \geq C_o$ on $T_o := [A-\tau, B+\tau]$.

In Case 1a, $W(\Delta) \leq \|\Delta\|_\infty^{\mathbb{R}} = 1$, whereas $\sigma(\Delta)^2 \geq C_o \int_c^{c+\delta} \Delta(x)^2 \, dx = C_o \delta^2/3$. Hence, (A.10) is satisfied if $K' \geq (3/C_o)^{1/2}$.

For Cases 1b and 2, we start with a more general consideration. Let $h(x) := 1\{x \in Q\}(\alpha + \gamma x)$ for real numbers $\alpha, \gamma$ and a non-degenerate interval $Q$ containing some point in $(c, c+\delta)$. Let $Q \cap T_o$ have end-points $x_o < y_o$. Elementary considerations then reveal that

$$\sigma(h)^2 \geq C_o \int_{x_o}^{y_o} (\alpha + \gamma x)^2 \, dx \geq \frac{C_o}{4}(y_o - x_o)(\|h\|_\infty^{T_o})^2.$$

We now deduce an upper bound for $W(h)/\|h\|_\infty^{T_o}$. If $Q \subset T_o$ or $\gamma = 0$, then $W(h)/\|h\|_\infty^{T_o} \leq 1$. Now, suppose that $\gamma \neq 0$ and $Q \not\subset T_o$. Then, $x_o, y_o \in T_o$ satisfy $y_o - x_o \geq \tau$ and, without loss of generality, let $\gamma = -1$. Now,

$$\|h\|_\infty^{T_o} = \max(|\alpha - x_o|, |\alpha - y_o|)$$
$$= (y_o - x_o)/2 + |\alpha - (x_o + y_o)/2|$$
$$\geq \tau/2 + |\alpha - (x_o + y_o)/2|.$$

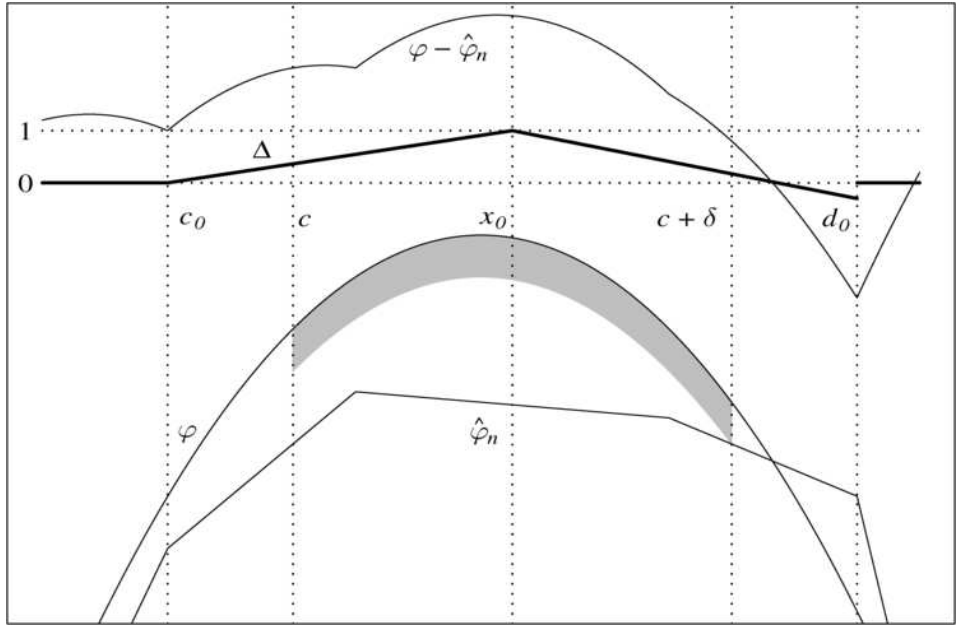

**Figure 5.** The perturbation function $\Delta$ in Case 2.

On the other hand, since $\varphi(x) \leq a_o - b_o|x|$ for certain constants $a_o, b_o > 0$,

$$\begin{aligned}
W(h) &\leq \sup_{x \in \mathbb{R}} \frac{|\alpha - x|}{\max(1, b_o|x| - a_o)} \\
&\leq \sup_{x \in \mathbb{R}} \frac{|\alpha| + |x|}{\max(1, b_o|x| - a_o)} \\
&= |\alpha| + (a_o + 1)/b_o \\
&\leq |\alpha - (x_o + y_o)/2| + (|A| + |B| + \tau)/2 + (a_o + 1)/b_o.
\end{aligned}$$

This entails that

$$\frac{W(h)}{\|h\|_\infty^{T_o}} \leq C_* := \frac{(|A| + |B| + \tau)/2 + (a_o + 1)/b_o}{\tau/2}.$$

In Case 1b, our function $\Delta$ is of the same type as $h$ above and $y_o - x_o \geq \delta$. Thus,

$$W(\Delta) \leq C_* \|h\|_\infty^{T_o} \leq 2 C_* C_o^{-1/2} \delta^{-1/2} \sigma(\Delta).$$

*Estimating log-concave densities*61

In Case 2, $\Delta$ may be written as $h_1 + h_2$, with two functions $h_1$ and $h_2$ of the same type as $h$ above having disjoint support and both satisfying $y_o - x_o \geq \delta/2$. Thus,

$$\begin{aligned} W(\Delta) &= \max(W(h_1), W(h_2)) \\ &\leq 2^{3/2} C_* C_o^{-1/2} \delta^{-1/2} \max(\sigma(h_1), \sigma(h_2)) \\ &\leq 2^{3/2} C_* C_o^{-1/2} \delta^{-1/2} \sigma(\Delta). \end{aligned}$$ □

To prove Lemma A.6, we need a simple exponential inequality.

**Lemma A.7.** *Let $Y$ be a random variable such that $\mathbb{E}(Y) = 0$, $\mathbb{E}(Y^2) = \sigma^2$ and $C := \mathbb{E}\exp(|Y|) < \infty$. Then, for arbitrary $t \in \mathbb{R}$,*

$$\mathbb{E}\exp(tY) \leq 1 + \frac{\sigma^2 t^2}{2} + \frac{C|t|^3}{(1 - |t|)_+}.$$

**Proof.**

$$\mathbb{E}\exp(tY) = \sum_{k=0}^{\infty} \frac{t^k}{k!} \mathbb{E}(Y^k) \leq 1 + \frac{\sigma^2 t^2}{2} + \sum_{k=3}^{\infty} \frac{|t|^k}{k!} \mathbb{E}(|Y|^k).$$

For any $y \geq 0$ and integers $k \geq 3$, $y^k e^{-y} \leq k^k e^{-k}$. Thus, $\mathbb{E}(|Y|^k) \leq \mathbb{E}\exp(|Y|) k^k e^{-k} = C k^k e^{-k}$. Since $k^k e^{-k} \leq k!$, which can be verified easily via induction on $k$,

$$\sum_{k=3}^{\infty} \frac{|t|^k}{k!} \mathbb{E}(|Y|^k) \leq C \sum_{k=3}^{\infty} |t|^k = \frac{C|t|^3}{(1 - |t|)_+}.$$ □

Lemma A.7 entails the following result for finite families of functions.

**Lemma A.8.** *Let $\mathcal{H}_n$ be a finite family of functions $h$ with $0 < W(h) < \infty$ such that $\#\mathcal{H}_n = O(n^p)$ for some $p > 0$. Then, for sufficiently large $D$,*

$$\lim_{n \to \infty} \mathbb{P}\left( \max_{h \in \mathcal{H}_n} \frac{|\int h\, d(\mathbb{F}_n - F)|}{\sigma(h)\rho_n^{1/2} + W(h)\rho_n^{2/3}} \geq D \right) = 0.$$

**Proof.** Since $W(ch) = cW(h)$ and $\sigma(ch) = c\sigma(h)$ for any $h \in \mathcal{H}_n$ and arbitrary constants $c > 0$, we may assume, without loss of generality, that $W(h) = 1$ for all $h \in \mathcal{H}_n$. Let $X$ be a random variable with log-density $\varphi$. Since

$$\limsup_{|x| \to \infty} \frac{\varphi(x)}{|x|} < 0$$

by Lemma A.1, the expectation of $\exp(t_o w(X))$ is finite for any fixed $t_o \in (0, 1)$, where $w(x) := \max(1, |\varphi(x)|)$. Hence,

$$\mathbb{E}\exp(t_o|h(X) - \mathbb{E}h(X)|) \leq C_o := \exp(t_o \mathbb{E}w(X))\mathbb{E}\exp(t_o w(X)) < \infty.$$

Lemma A.7, applied to $Y := t_o(h(X) - \mathbb{E}h(X))$, implies that

$$\mathbb{E}\exp[t(h(X) - \mathbb{E}h(X))] = \mathbb{E}((t/t_o)Y) \le 1 + \frac{\sigma(h)^2 t^2}{2} + \frac{C_1|t|^3}{(1-C_2|t|)_+}$$

for arbitrary $h \in \mathcal{H}_n$, $t \in \mathbb{R}$ and constants $C_1, C_2$ depending on $t_o$ and $C_o$. Consequently,

$$\begin{aligned}
\mathbb{E}\exp\left(t\int h\,d(\mathbb{F}_n - F)\right) &= \mathbb{E}\exp\left((t/n)\sum_{i=1}^n (h(X_i) - \mathbb{E}h(X))\right) \\
&= (\mathbb{E}\exp((t/n)(h(X) - \mathbb{E}h(X))))^n \\
&\le \left(1 + \frac{\sigma(h)^2 t^2}{2n^2} + \frac{C_1|t|^3}{n^3(1 - C_2|t|/n)_+}\right)^n \\
&\le \exp\left(\frac{\sigma(h)^2 t^2}{2n} + \frac{C_1|t|^3}{n^2(1-C_2|t|/n)_+}\right).
\end{aligned}$$

It now follows from Markov's inequality that

$$\mathbb{P}\left(\left|\int h\,d(\mathbb{F}_n - F)\right| \ge \eta\right) \le 2\exp\left(\frac{\sigma(h)^2 t^2}{2n} + \frac{C_1 t^3}{n^2(1-C_2 t/n)_+} - t\eta\right) \qquad (A.12)$$

for arbitrary $t, \eta > 0$. Specifically, let $\eta = D(\sigma(h)\rho_n^{1/2} + \rho_n^{2/3})$ and set

$$t := \frac{n\rho_n^{1/2}}{\sigma(h) + \rho_n^{1/6}} \le n\rho_n^{1/3} = o(n).$$

Then, the bound (A.12) is not greater than

$$2\exp\left(\frac{\sigma(h)^2 \log n}{2(\sigma(h) + \rho_n^{1/6})^2} + \frac{C_1 \rho_n^{1/2} \log n}{(\sigma(h) + \rho_n^{1/6})^3(1 - C_2\rho_n^{1/3})_+} - D\log n\right)$$

$$\le 2\exp\left[\left(\frac{1}{2} + \frac{C_1}{(1-C_2\rho_n^{1/3})_+} - D\right)\log n\right] = 2\exp((O(1) - D)\log n).$$

Consequently, for sufficiently large $D > 0$,

$$\mathbb{P}\left(\max_{h\in\mathcal{H}_n} \frac{|\int h\,d(\mathbb{F}_n - F)|}{\sigma(h)\rho_n^{1/2} + W(h)\rho_n^{2/3}} \ge D\right)$$

$$\le \#\mathcal{H}_n 2\exp((O(1) - D)\log n) = O(1)\exp((O(1) + p - D)\log n) \to 0. \qquad \square$$

**Proof of Lemma A.6.** Let $\mathcal{H}$ be the family of all functions $h$ of the form

$$h(x) = 1\{x \in Q\}(c + dx),$$

with any interval $Q \subset \mathbb{R}$ and real constants $c, d$ such that $h$ is non-negative. Suppose that there exists a constant $C = C(f)$ such that

$$\mathbb{P}\left(\sup_{h \in \mathcal{H}} \frac{|\int h \, d(\mathbb{F}_n - F)|}{\sigma(h)\rho_n^{1/2} + W(h)\rho_n^{2/3}} \leq C\right) \to 1. \tag{A.13}$$

For any $m \in \mathbb{N}$, an arbitrary function $\Delta \in \mathcal{D}_m$ may be written as

$$\Delta = \sum_{i=1}^{M} h_i$$

with $M = 2m + 2$ functions $h_i \in \mathcal{H}$ having pairwise disjoint supports. Consequently,

$$\sigma(\Delta) = \left(\sum_{i=1}^{M} \sigma(h_i)^2\right)^{1/2} \geq M^{-1/2} \sum_{i=1}^{M} \sigma(h_i),$$

by the Cauchy–Schwarz inequality, while

$$W(\Delta) = \max_{i=1,\ldots,M} W(h_i) \geq M^{-1} \sum_{i=1}^{M} W(h_i).$$

Consequently, (A.13) entails that

$$\left|\int \Delta \, d(\mathbb{F}_n - F)\right| \leq \sum_{i=1}^{M} \left|\int h_i \, d(\mathbb{F}_n - F)\right|$$

$$\leq C\left(\sum_{i=1}^{M} \sigma(h_i)\rho_n^{1/2} + \sum_{i=1}^{M} W(h_i)\rho_n^{2/3}\right)$$

$$\leq 4C(\sigma(\Delta)m^{1/2}\rho_n^{1/2} + W(\Delta)m\rho_n^{2/3})$$

uniformly in $m \in \mathbb{N}$ and $\Delta \in \mathcal{D}_m$, with probability tending to one as $n \to \infty$.

It remains to verify (A.13). To this end, we use a bracketing argument. With the weight function $w(x) = \max(1, |\varphi(x)|)$, let $-\infty = t_{n,0} < t_{n,1} < \cdots < t_{n,N(n)} = \infty$ such that for $I_{n,j} := (t_{n,j-1}, t_{n,j}]$,

$$(2n)^{-1} \leq \int_{I_{n,j}} w(x)^2 f(x) \, dx \leq n^{-1} \quad \text{for } 1 \leq j \leq N(n),$$

with equality if $j < N(n)$. Since $1 \leq \int \exp(t_o w(x)) f(x) \, dx < \infty$, such a partition exists with $N(n) = O(n)$. For any $h \in \mathcal{H}$, we define functions $h_{n,\ell}, h_{n,u}$ as follows. Let $\{j, \ldots, k\}$ be the set of all indices $i \in \{1, \ldots, N(n)\}$ such that $\{h > 0\} \cap I_{n,i} \neq \emptyset$. We then define

$$h_{n,\ell}(x) := 1_{\{t_{n,j} < x \leq t_{n,k-1}\}} h(x)$$

and

$$h_{n,u}(x) := h_{n,\ell}(x) + 1\{x \in I_{n,j} \cup I_{n,k}\}W(h)w(x).$$

Note that $0 \leq h_{n,\ell} \leq h \leq h_{n,u} \leq W(h)w$. Consequently, $W(h_{n,\ell}) \leq W(h) = W(h_{n,u})$. Suppose, for the moment, that the assertion is true for the (still infinite) family $\mathcal{H}_n := \{h_{n,\ell}, h_{n,u} : h \in \mathcal{H}\}$ in place of $\mathcal{H}$. It then follows from $w \geq 1$ that

$$\begin{aligned}
\int h \, \mathrm{d}(\mathbb{F}_n - F) &\leq \int h_{n,u} \, \mathrm{d}\mathbb{F}_n - \int h_{n,\ell} \, \mathrm{d}F \\
&= \int h_{n,u} \, \mathrm{d}(\mathbb{F}_n - F) + \int (h_{n,u} - h_{n,\ell}) \, \mathrm{d}F \\
&\leq \int h_{n,u} \, \mathrm{d}(\mathbb{F}_n - F) + W(h) \int_{I_{n,j} \cup I_{n,k}} w(x)^2 \, \mathrm{d}F \\
&\leq \int h_{n,u} \, \mathrm{d}(\mathbb{F}_n - F) + 2W(h)n^{-1} \\
&\leq C(\sigma(h_{n,u})\rho_n^{1/2} + \rho_n^{2/3}) + 2n^{-1} \\
&\leq C(\sigma(h)\rho_n^{1/2} + 2^{1/2}W(h)n^{-1/2}\rho_n^{1/2} + \rho_n^{2/3}) + 2W(h)n^{-1} \\
&\leq (C + o(1))(\sigma(h)\rho_n^{1/2} + W(h)\rho_n^{2/3}),
\end{aligned}$$

uniformly in $h \in \mathcal{H}$ with asymptotic probability one. Analogously,

$$\begin{aligned}
\int h \, \mathrm{d}(\mathbb{F}_n - F) &\geq \int h_{n,\ell} \, \mathrm{d}(\mathbb{F}_n - F) - 2W(h)n^{-1} \\
&\geq -C(\sigma(h_{n,\ell})\rho_n^{1/2} + W(h)\rho_n^{2/3}) - 2W(h)n^{-1} \\
&\geq -(C + o(1))(\sigma(h)\rho_n^{1/2} + W(h)\rho_n^{2/3}),
\end{aligned}$$

uniformly in $h \in \mathcal{H}$ with asymptotic probability one.

To accord with Lemma A.8, we must now deal with $\mathcal{H}_n$. For any $h \in \mathcal{H}$, the function $h_{n,\ell}$ may be written as

$$h(t_{n,j})g^{(1)}_{n,j,k} + h(t_{n,k-1})g^{(2)}_{n,j,k},$$

with the "triangular functions"

$$g^{(1)}_{n,j,k}(x) := \frac{t_{n,k-1} - x}{t_{n,k-1} - t_{n,j}}$$

and

$$g^{(2)}_{n,j,k}(x) := \frac{x - t_{n,j}}{t_{n,k-1} - t_{n,j}} \qquad \text{for } 1 \leq j < k \leq N(n), k - j \geq 2.$$

In case of $k - j \leq 1$, we set $g_{n,j,k}^{(1)} := g_{n,j,k}^{(2)} := 0$. Moreover,

$$h_{n,u} = h_{n,\ell} + W(h)g_{n,j} + 1\{k > j\}W(h)g_{n,k},$$

with $g_{n,i}(x) := 1\{x \in I_{n,i}\}w(x)$. Consequently, all functions in $\mathcal{H}_n$ are linear combinations with *non-negative* coefficients of at most four functions in the finite family

$$\mathcal{G}_n := \{g_{n,i} : 1 \leq i \leq N(n)\} \cup \{g_{n,j,k}^{(1)}, g_{n,j,k}^{(2)} : 1 \leq j < k \leq N(n)\}.$$

Since $\mathcal{G}_n$ contains $O(n^2)$ functions, it follows from Lemma A.8 that for some constant $D > 0$,

$$\left| \int g \, \mathrm{d}(\mathbb{F}_n - F) \right| \leq D(\sigma(g)\rho_n^{1/2} + W(g)\rho_n^{2/3})$$

for all $g \in \mathcal{G}_n$ with asymptotic probability one. The assertion about $\mathcal{H}_n$ now follows from the basic observation that for $h = \sum_{i=1}^{4} \alpha_i g_i$ with non-negative functions $g_i$ and coefficients $\alpha_i \geq 0$,

$$\sigma(h) \geq \left( \sum_{i=1}^{4} \alpha_i^2 \sigma(g_i)^2 \right)^{1/2} \geq 2^{-1} \sum_{i=1}^{4} \alpha_i \sigma(g_i),$$

$$W(h) \geq \max_{i=1,\ldots,4} \alpha_i W(g_i) \geq 4^{-1} \sum_{i=1}^{4} \alpha_i W(g_i). \qquad \square$$

## A.4. Proofs for the gap problem and of $\hat{F}_n$'s consistency

**Proof of Theorem 4.3.** Suppose that $\hat{\varphi}_n$ is linear on an interval $[a, b]$. Then, for $x \in [a, b]$ and $\lambda_x := (x - a)/(b - a) \in [0, 1]$,

$$\varphi(x) - (1 - \lambda_x)\varphi(a) - \lambda_x\varphi(b)$$
$$= (1 - \lambda_x)(\varphi(x) - \varphi(a)) - \lambda_x(\varphi(b) - \varphi(x))$$
$$= (1 - \lambda_x) \int_a^x \varphi'(t) \, \mathrm{d}t - \lambda_x \int_x^b \varphi'(t) \, \mathrm{d}t$$
$$= (1 - \lambda_x) \int_a^x (\varphi'(t) - \varphi'(x)) \, \mathrm{d}t + \lambda_x \int_x^b (\varphi'(x) - \varphi'(t)) \, \mathrm{d}t$$
$$\geq C(1 - \lambda_x) \int_a^x (x - t) \, \mathrm{d}t + C\lambda_x \int_x^b (t - x) \, \mathrm{d}t$$
$$= C(b - a)^2 \lambda_x(1 - \lambda_x)/2$$
$$= C(b - a)^2/8 \qquad \text{if } x = x_o := (a + b)/2.$$

This entails that $\sup_{[a,b]} |\hat\varphi_n - \varphi| \geq C(b-a)^2/16$. For if $\hat\varphi_n < \varphi + C(b-a)^2/16$ on $\{a,b\}$, then

$$\begin{aligned}
\varphi(x_o) - \hat\varphi_n(x_o) &= \varphi(x_o) - (\hat\varphi_n(a) + \hat\varphi_n(b))/2 \\
&> \varphi(x_o) - (\varphi(a) + \varphi(b))/2 - C(b-a)^2/16 \\
&\geq C(b-a)^2/8 - C(b-a)^2/16 = C(b-a)^2/16.
\end{aligned}$$

Consequently, if $|\hat\varphi_n - \varphi| \leq D_n \rho_n^{\beta/(2\beta+1)}$ on $T_n := [A + \rho_n^{1/(2\beta+1)}, B - \rho_n^{1/(2\beta+1)}]$ with $D_n = O_{\mathrm{p}}(1)$, then the longest subinterval of $T_n$ containing no points from $\mathcal{S}_n$ has length at most $4 D_n^{1/2} C^{-1/2} \rho_n^{\beta/(4\beta+2)}$. Since $T_n$ and $T = [A, B]$ differ by two intervals of length $\rho_n^{1/(2\beta+1)} = O(\rho_n^{\beta/(4\beta+2)})$, these considerations yield the assertion about $\mathcal{S}_n(\hat\varphi_n)$.  $\square$

**Proof of Theorem 4.4.** Let $\delta_n := \rho_n^{1/(2\beta+1)}$ and $r_n := D\rho_n^{\beta/(4\beta+2)} = D\delta_n^{1/2}$ for some constant $D > 0$. Since $r_n \to 0$ but $nr_n \to \infty$, it follows from boundedness of $f$ and a theorem of Stute (1982) about the modulus of continuity of univariate empirical processes that

$$\begin{aligned}
\omega_n &:= \sup_{x,y \in \mathbb{R}: |x-y| \leq r_n} |(\mathbb{F}_n - F)(x) - (\mathbb{F}_n - F)(y)| \\
&= O_{\mathrm{p}}(n^{-1/2} r_n^{1/2} \log(1/r_n)^{1/2}) \\
&= O_{\mathrm{p}}(\rho_n^{(5\beta+2)/(8\beta+4)}).
\end{aligned}$$

If $D$ is sufficiently large, the asymptotic probability that for any point $x \in [A+\delta_n, B-\delta_n]$, there exists a point $y \in \mathcal{S}_n(\hat\varphi_n) \cap [A+\delta_n, B-\delta_n]$ with $|x-y| \leq r_n$, is equal to one. In that case, it follows from Corollary 2.5 and Theorem 4.1 that

$$\begin{aligned}
|(\hat F_n - \mathbb{F}_n)(x)| &\leq |(\hat F_n - \mathbb{F}_n)(x) - (\hat F_n - \mathbb{F}_n)(y)| + n^{-1} \\
&\leq |(\hat F_n - F)(x) - (\hat F_n - F)(y)| + \omega_n + n^{-1} \\
&\leq \int_{\min(x,y)}^{\max(x,y)} |\hat f_n - f|(x)\,\mathrm{d}x + \omega_n + n^{-1} \\
&\leq O_{\mathrm{p}}(r_n \rho_n^{\beta/(2\beta+1)}) + \omega_n + n^{-1} \\
&= O_{\mathrm{p}}(\rho_n^{3\beta/(4\beta+2)}). \qquad \square
\end{aligned}$$

## Acknowledgements

This work is part of the second author's PhD dissertation, written at the University of Bern. The authors thank an anonymous referee for valuable remarks and some important references. This work was supported by the Swiss National Science Foundation.